# Improved Family Of Estimators Of Population Mean In Simple Random Sampling

Rajesh Singh, Mukesh Kumar and Manoj K. Chaudhary

Department of Statistics, Banaras Hindu University (U.P.), India

rsinghstat@yahoo.co.in

**Abstract**

In this paper, a procedure is given for estimating the population mean in simple random sampling without replacement in the presence of auxiliary information. The mean squared error expressions of the proposed estimators have been derived under large sample approximation. We have compared the performances of the proposed estimators with several existing estimators. Both theoretical and empirical findings are encouraging and support the soundness of the proposed procedure for mean estimation.

1. Introduction

The problem of estimating the population mean in the presence of an auxiliary variable has been widely discussed in finite population sampling literature. Ratio, product and difference methods of estimation are good examples in this context. Ratio method of estimation is quite effective when there is high positive correlation between study and auxiliary variables. On the other hand, if correlation is negative (high), the product method of estimation can be employed efficiently.

In recent years, a number of research papers on ratio-type, exponential ratio-type and regression-type estimators have appeared, based on different types of transformations. Some important contributions in this area are due to Singh and Tailor [12], Shabbir and Gupta [6], Kadilar and Cingi [3,4], Khosnevisan et al. [5], Singh et al. [7], Singh et al. [8] and Singh et al. [10]. We begin by introducing some terminology used in the article.

Let $\Omega$ be a finite population consisting of N units from which a sample of size n is to be drawn by simple random sampling without replacement (SRSWOR). Let y and x denote the study and auxiliary variables having sample means $\bar{y}$ and $\bar{x}$ respectively. It is desired to estimate the population mean $\bar{Y}$ of the study variate y using information on the population mean $\bar{X}$ of the auxiliary variate x.

Let $e_0 = \frac{(\bar{y}-\bar{Y})}{\bar{Y}}$ and $e_1 = \frac{(\bar{x}-\bar{X})}{\bar{X}}$.

This gives E ($e_0$) = E ($e_1$) = 0,

$$E(e_0^2) = \frac{(1-f)}{n} C_y^2, \quad E(e_1^2) = \frac{(1-f)}{n} C_x^2 \quad \text{and} \quad E(e_0 e_1) = \left(\frac{1-f}{n}\right) \rho C_y C_x,$$

where $f = \frac{n}{N}$, $C_y = \frac{S_y}{\bar{Y}}$, $C_x = \frac{S_x}{\bar{X}}$ and $\rho = \frac{S_{xy}}{S_x S_y}$.

We discuss below some of the estimators available in the literature.

The most common estimator of $\bar{Y}$ is the sample mean estimator defined as

$$\bar{y}_0 = \bar{y}. \tag{1}$$

The mean square error (MSE) of $\bar{y}_0$ is given by

$$MSE(\bar{y}_0) = \frac{1-f}{n} \bar{Y}^2 C_y^2 \tag{2}$$

The usual ratio estimator for $\bar{Y}$ is

$$\bar{y}_R = \bar{y}\left(\frac{\bar{X}}{\bar{x}}\right) \tag{3}$$

The MSE of $\bar{y}_R$, to first –order of approximation is given by

$$MSE(\bar{y}_R) = \frac{1-f}{n} \bar{Y}^2 \left[C_y^2 + C_x^2 - 2\rho C_y C_x\right] \tag{4}$$

## 2. Estimators proposed by other authors

Singh and Tailor [9] estimator

By using the known correlation coefficient ($\rho$), Singh and Tailor [9] introduced the following ratio type estimator

$$\bar{y}_{ST} = \bar{y}\left(\frac{\bar{X}+\rho}{\bar{x}+\rho}\right) \tag{5}$$

with MSE

$$\text{MSE}(\bar{y}_{ST}) = \frac{1-f}{n} \bar{Y}^2 [C_y^2 + \eta^2 C_x^2 - 2\eta\rho C_y C_x] \qquad (6)$$

where

$$\eta = \frac{\bar{X}}{(\bar{X}+\rho)},$$

Khosnevisan et al. [5] suggested following modified estimator

$$\bar{y}_{KH} = \bar{y}\left[\frac{a\bar{X}+b}{\alpha(a\bar{x}+b)+(1-\alpha)(a\bar{X}+b)}\right]^g \qquad (7)$$

where g and $\alpha$ are suitable constants, a and b are either real numbers or the functions of the known parameters of the auxiliary variable, x, such as coefficient of variation ($C_x$), kurtosis ($\beta_2(x)$) and correlation coefficient ($\rho$).

**Remark 1**- Here we would like to mention that the choice of the estimator depends on the availability and values of the various parameter(s) used (for choice of the parameters a and b refer to Khosnevisan et al. [5] and Singh et al. [8]).

The minimum MSE of $\bar{y}_{KH}$, to first order of approximation, is given by

$$\text{MSE}(\bar{y}_{KH})_{min} = \frac{1-f}{n}\bar{Y}^2 C_y^2(1-\rho^2) \qquad (8)$$

Bahl and Tuteja [1] suggested an exponential ratio- type estimator defined as

$$\bar{y}_{BT} = \bar{y}\exp\left[\frac{\bar{X}-\bar{x}}{\bar{X}-\bar{x}}\right] \qquad (9)$$

The MSE of $\bar{y}_{BT}$, to the first order of approximation, is given by

$$\text{MSE}(\bar{y}_{BT}) = \frac{1-f}{n}\bar{Y}^2\left(C_y^2 + \frac{C_x^2}{4} - \rho C_y C_x\right) \qquad (10)$$

Singh et al. [7] suggested a modified form of Bahl and Tuteja [1] estimator

$$\bar{y}_8 = \bar{y} \exp\left[\frac{(a\bar{X}+b)-(a\bar{x}+b)}{(a\bar{X}+b)+(a\bar{x}+b)}\right] \tag{11}$$

The minimum MSE of $\bar{y}_8$ to the first order of approximation is given by

$$\text{MSE}(\bar{y}_8)_{\min} = \frac{1-f}{n}\bar{Y}^2 C_y^2 (1-\rho^2) \tag{12}$$

Note that MSE in (7) and (11) is equal to the MSE of the linear regression estimator

$$\bar{y}_{reg} = \bar{y} + b(\bar{X}-\bar{x}).$$

Kadilar and Cingi [3] introduced following class of estimators for estimating $\bar{Y}$ -

$$\bar{y}_{KCi} = [\bar{y}+b(\bar{X}-\bar{x})]\gamma_i, \quad i=1\ldots5, \tag{13}$$

where $\gamma_i$'s are defined as

$$\gamma_1 = \frac{\bar{X}}{\bar{x}}, \qquad \gamma_2 = \frac{\bar{X}+C_x}{\bar{x}+C_x}, \qquad \gamma_3 = \frac{\bar{X}+\beta_2(x)}{\bar{x}+\beta_2(x)},$$

$$\gamma_4 = \frac{\bar{X}\beta_2(x)+C_x}{\bar{x}\beta_2(x)+C_x}, \qquad \gamma_5 = \frac{\bar{X}C_x+\beta_2(x)}{\bar{x}C_x+\beta_2(x)}.$$

The MSE of $\bar{y}_{KCi}$ (i=1, 2..5), to first order of approximation, is given by

$$\text{MSE}(\bar{y}_{KCi}) = \frac{1-f}{n}\bar{Y}^2[\gamma_i^{*2} C_x^2 + C_y^2(1-\rho^2)] \tag{14}$$

where $\gamma_1^* = 1,$ $\gamma_2^* = \bar{X}/(\bar{X}+C_x),$ $\gamma_3^* = \bar{X}/(\bar{X}+\beta_2(x)),$

$\gamma_4^* = \bar{X}\beta_2(x)/(\bar{X}\beta_2(x)+C_x),$ and $\gamma_5^* = \bar{X}C_x/(\bar{X}+\beta_2(x)).$

Using Singh and Tailor [9] transformation, Kadilar and Cingi [4] introduced another class of estimators, given by

$$\bar{y}_{KCi}^* = [\bar{y} + b(\bar{X} - \bar{x})]\psi_i, \quad i=1,2,\ldots,5. \tag{15}$$

$$\psi_1 = \frac{\bar{X}+\rho}{\bar{x}+\rho}, \quad \psi_2 = \frac{\bar{X}C_x+\rho}{\bar{x}C_x+\rho}, \quad \psi_3 = \frac{\bar{X}\rho+C_x}{\bar{x}\rho+C_x},$$

$$\psi_4 = \frac{\bar{X}\beta_2(x)+\rho}{\bar{x}\beta_2(x)+\rho}, \quad \psi_5 = \frac{\bar{X}\rho+\beta_2(x)}{\bar{x}\rho+\beta_2(x)}.$$

The MSE of $\bar{y}_{KCi}^*$, to first order of approximation, is given by

$$MSE(\bar{y}_{KCi}^*) = \frac{1-f}{n}\bar{Y}^2\left[\psi_i^{*2}C_x^2 + C_y^2(1-\rho^2)\right], \tag{16}$$

where

$$\psi_1^* = \frac{\bar{X}}{(\bar{X}+\rho)}, \quad \psi_2^* = \bar{X}C_x/(\bar{X}C_x+\rho),$$

$$\psi_3^* = \frac{\bar{X}\rho}{(\bar{X}\rho+C_x)}, \quad \psi_4^* = \bar{X}\beta_2(x)/(\bar{X}\beta_2(x)+\rho),$$

and $\psi_5^* = \frac{\bar{X}\rho}{(\bar{X}\rho+\beta_2(x))}$.

Gupta and Shabbir [2] introduced following general class of ratio-type estimators-

$$\bar{y}_{GS} = [\omega_1\bar{y} + \omega_2(\bar{X} - \bar{x})]\left(\frac{a\bar{X}+b}{a\bar{x}+b}\right) \tag{17}$$

where $\omega_1$ and $\omega_2$ are weights.

The minimum MSE of $\bar{y}_{GS}$, to first order of approximation is given by

$$\text{MSE}(\bar{y}_{GS})_{min} = \frac{\text{MSE}(\bar{y}_{reg})}{1 + \left(\frac{1-f}{n}\right)C_y^2(1-\rho^2)} \qquad (18)$$

3. **Proposed estimator**

Motivated by Gupta and Shabbir [2], we propose the following ratio type estimator-

$$\bar{y}_\lambda = [q_1\bar{y} + q_2(\bar{X} - \bar{x})]\left[\lambda\left(\frac{a\bar{X}+b}{a\bar{x}+b}\right) + (1-\lambda)\exp\left\{\frac{a\bar{X}-a\bar{x}}{a\bar{X}+a\bar{x}+2b}\right\}\right] \qquad (19)$$

where $q_1$ and $q_2$ are weights whose values are to be determined later, $\lambda$ is a suitably chosen constant and a and b are same as defined in (7).

**Remark 2** – For the choice $q_1 = 1, q_2 = 0$, a=1, b=0 and $\lambda = 1$ the estimator $\bar{y}_\lambda$ changes to the usual estimator $\bar{y}_R$. Again for the choice $q_1 = 1, q_2 = 0$, a=1, b=0 and $\lambda = 0$ the estimator $\bar{y}_\lambda$ changes to Bahl and Tuteja [1] estimator $\bar{y}_{BT}$. On selecting $\lambda = 1$, the estimator $\bar{y}_\lambda$ reduces to the Gupta and Shabbir estimator [2] $\bar{y}_{GS}$. Several other estimators can be generated from the estimator $\bar{y}_\lambda$ by suitable choice of $q_1, q_2$, a, b and $\lambda$.

From (19), we have

$$\bar{y}_\lambda = [q_1\bar{y}(1+e_0) - q_2(\bar{X}e_1)]\left[\lambda(1+\theta e_1)^{-1} + (1-\lambda)\exp\left\{-\frac{\theta}{2}e_1\left(1+\frac{\theta e_1}{2}\right)^{-1}\right\}\right]$$

$$= [q_1\bar{y}(1+e_0) - q_2(\bar{X}e_1)]\left[\left(1 - \frac{\theta e_1}{2} + \frac{3}{8}\theta^2 e_1^2\right) + \lambda\left(\frac{-\theta e_1}{2} + \frac{5}{8}\theta^2 e_1^2\right)\right]$$

$$= q_1\bar{y}\left\{1 + e_0 - \frac{\theta e_1}{2}(1+\lambda) + \frac{1}{8}\theta^2 e_1^2(3+5\lambda) - \frac{\theta}{2}e_0 e_1(1+\lambda)\right\}$$

$$\qquad - q_2\bar{X}\left\{e_1 - \frac{\theta}{2}e_1^2 + \lambda\left(\frac{-\theta e_1^2}{2}\right)\right\} \qquad (20)$$

where $\theta = \frac{a\bar{X}}{a\bar{X}+b}$

From (20), we have

$$(\bar{y}_\lambda - \bar{Y}) = (q_1 - 1)\bar{Y} + q_1\bar{Y}\left\{e_0 - \frac{\theta e_1(1+\lambda)}{2} + \frac{\theta^2 e_1^2}{8}(3+5\lambda) - \frac{\theta}{2}e_0 e_1(1+\lambda)\right\}$$

$$-q_2\bar{X}\left\{e_1 - \frac{\theta e_1^2}{2} - \frac{\theta \lambda e_1^2}{2}\right\} \qquad (21)$$

Using (21), we get the MSE of $\bar{y}_\lambda$, to first order of approximation, as given by

$$MSE(\bar{y}_\lambda) = E(\bar{y}_\lambda - \bar{Y})^2$$

$$\cong E\left[(q_1 - 1)\bar{Y} + q_1\bar{Y}\left\{e_0 - \frac{\theta e_1}{2}(1+\lambda)\right\} - q_2\bar{X}e_1\right]^2$$

$$\cong E\left[(q_1 - 1)\bar{Y} + q_1\bar{Y}e_0 - \left(\frac{\theta}{2}q_1\bar{Y}(1+\lambda) + q_2\bar{X}\right)e_1\right]^2$$

$$\cong ((q_1 - 1)\bar{Y})^2 + (q_1\bar{Y})^2 E(e_0^2) + \left(\frac{\theta}{2}q_1\bar{Y}(1+\lambda) + q_2\bar{X}\right)^2 E(e_1^2)$$

$$-2q_1\bar{Y}\left(\frac{\theta}{2}q_1\bar{Y}(1+\lambda) + q_2\bar{X}\right)E(e_1 e_0)$$

$$\cong ((q_1 - 1)\bar{Y})^2 + (q_1\bar{Y})^2 \frac{1-f}{n}C_y^2 + \left(\frac{\theta}{2}q_1\bar{Y}(1+\lambda) + q_2\bar{X}\right)^2 \frac{1-f}{n}C_x^2$$

$$-2q_1\bar{Y}\left(\frac{\theta}{2}q_1\bar{Y}(1+\lambda) + q_2\bar{X}\right)\frac{1-f}{n}\rho C_y C_x$$

On rearranging the terms we get

$$MSE(\bar{y}_\lambda) \cong \bar{Y}^2(q_1 - 1)^2 + \frac{1-f}{n}\left[q_1^2\bar{Y}^2\left\{C_y^2 + \frac{\theta^2 C_x^2}{4}(1+\lambda)^2 - (1+\lambda)\rho C_y C_x\right\} + q_2^2\bar{X}^2 C_x^2 - 2q_1 q_2 \bar{Y}\bar{X}C_x\left\{\rho C_y - \frac{(1+\lambda)}{2}\theta C_x\right\}\right]$$

(22)

Minimization of (22) with respect to $q_i (i = 1,2)$ yields optimum values of $q_1$ and $q_2$ as

$$q_{10} = \frac{1}{1+\left(\frac{1-f}{n}\right)C_y^2(1-\rho^2)}$$

and

$$q_{20} = q_{10}\cdot\left(\frac{Y}{XC_x}\right)\left\{\rho C_y - \left(\frac{1+\lambda}{2}\right)\theta C_x\right\} \tag{23}$$

We assume that the unknown parameters $C_y$ and $\rho$ are easily estimable from the preliminary data as in Tracy et al. [12] and Singh et al. [8]. Substituting the optimum values $q_{10}$ and $q_{20}$ in (22), the minimum MSE of $\bar{y}_p$, to first order of approximation, is given by

$$MSE(\bar{y}_\lambda)_{min} = \frac{MSE(\bar{y}_{reg})}{1+\left(\frac{1-f}{n}\right)C_y^2(1-\rho^2)} \tag{24}$$

## 4. Efficiency comparisons

First we compare the efficiency of proposed estimator with mean per unit estimator

$$MSE(\bar{y}_0) - MSE(\bar{y}_\lambda)_{min} \geq 0$$

$$f_1 Y^2 C_y^2 - \frac{f_1 Y^2 C_y^2(1-\rho^2)}{1+f_1 C_y^2(1-\rho^2)} \geq 0$$

$$f_1 C_y^2(1-\rho^2) + \rho^2 \geq 0 \tag{25}$$

where

$$f_1 = \left(\frac{1-f}{n}\right).$$

Next we compare the proposed estimator with usual ratio estimator

$$MSE(\bar{y}_R) - MSE(\bar{y}_\lambda)_{min} \geq 0$$

$$f_1 \bar{Y}^2(C_y^2 + C_x^2 - 2\rho C_y C_x) - \frac{f_1 \bar{Y}^2 C_y^2(1-\rho^2)}{1+f_1 C_y^2(1-\rho^2)} \geq 0$$

$$f_1 \bar{Y}^2\left[(C_y^2 + C_x^2 - 2\rho C_y C_x) - \frac{C_y^2(1-\rho^2)}{1+f_1 C_y^2(1-\rho^2)}\right] \geq 0 \tag{26}$$

Similarly, we can compare the proposed estimator with exponential ratio type estimator (with a=1 and b=0)-

$$MSE(\bar{y}_g) - MSE(\bar{y}_\lambda)_{min} \geq 0$$

$$f_1 \bar{Y}^2\left(C_y^2 + \frac{C_x^2}{4} - \rho C_y C_x\right) - \frac{f_1 \bar{Y}^2 C_y^2(1-\rho^2)}{1+f_1 C_y^2(1-\rho^2)} \geq 0$$

$$f_1 \bar{Y}^2\left[\left(C_y^2 + \frac{C_x^2}{4} - \rho C_y C_x\right) - \frac{C_y^2(1-\rho^2)}{1+f_1 C_y^2(1-\rho^2)}\right] \geq 0 \tag{27}$$

We see, conditions (26) and (27) are always satisfied.

Next we compare the efficiency of proposed estimator with Kadilar and Cingi [3] estimator $\bar{y}_{KC1}$.

$$MSE(\bar{y}_{KC1}) - MSE(\bar{y}_\lambda)_{min} \geq 0$$

$$f_1 \bar{Y}^2[\gamma_1^{*2} C_x^2 + C_y^2(1-\rho^2)] - \frac{f_1 \bar{Y}^2 C_y^2(1-\rho^2)}{1+f_1 C_y^2(1-\rho^2)} \geq 0$$

On solving, we get final expressions as

$$f_1 \bar{Y}^2 \gamma_1^{*2} C_x^2\{1 + f_1 C_y^2(1-\rho^2)\} + \{f_1 \bar{Y} C_y^2(1-\rho^2)\}^2 \geq 0 \tag{28}$$

This is always true.

Similarly we compare the proposed estimator with Kadilar and Cingi [4] estimator $\bar{y}_{KC1}^*$.

The only change in expression is to replace $\psi_i^{*2}$ by $\psi_i^{*2}$, we get

$$f_1 Y^2 \psi_i^{*2} C_x^2 (1 + f_1 C_y^2(1-\rho^2)) + \{f_1 Y C_y^2(1-\rho^2)\}^2 \geq 0 \qquad (29)$$

This is also always satisfied.

Next we compare proposed estimator with ( Khosnevisan et al.[5], Singh et al. [7] and regression estimators) as-

$$\left[f_1 Y^2 C_y^2(1-\rho^2) - \frac{f_1 Y^2 C_y^2(1-\rho^2)}{1 + f_1 C_y^2(1-\rho^2)}\right] \geq 0$$

$$\left[f_1 Y C_y^2(1-\rho^2)\right]^2 \geq 0. \qquad (30)$$

This is true always.

Finally we compare the efficiency of proposed estimator with Gupta and Shabbir [2] estimator.

$$\text{MSE}(\bar{y}_{GS}) - \text{MSE}(\bar{y}_\lambda)_{min} \geq 0$$

$$\frac{f_1 Y^2 C_y^2(1-\rho^2)}{1 + f_1 C_y^2(1-\rho^2)} - \frac{f_1 Y^2 C_y^2(1-\rho^2)}{1 + f_1 C_y^2(1-\rho^2)} = 0 \qquad (31)$$

as we see both the MSE(s) are equal.

**5. Empirical study.**

For empirical study we use two data sets earlier used by Singh and Tailor [9] and Singh and Agnihotri [11].

**Table 1:**

**Data statistics:**

| Population | N | n | $\bar{Y}$ | $\bar{X}$ | $C_y$ | $C_x$ | $\rho$ | $\beta_2(x)$ |
|---|---|---|---|---|---|---|---|---|
| Population 1 | 200 | 50 | 500 | 25 | 15 | 2 | 0.90 | 50 |
| Population 2 | 278 | 30 | 39.068 | 25.111 | 1.445 | 1.620 | 0.721 | 38.890 |

**Table 2: Percent relative efficiency (PRE) values of estimators**

| | | | | The value of $\lambda$ | | | |
|---|---|---|---|---|---|---|---|
| The values of | | | | 1 | | 0 | |
| $q_1$ | $q_2$ | Estimator | Values of a, b | PRE1 | PRE2 | PRE1 | PRE2 |
| 1 | 0 | $\bar{y}_R$ | a=1, b=0 | 100 | 100 | 100 | 100 |
| 1 | 0 | $\bar{y}_{CR}$ | $a = 1, b = 0$ | 128.571 | 156.190 | 113.065 | 197.667 |
| | | $\bar{y}_{SD}$ | $a = 1, b = C_x$ | 126.100 | 169.350 | 112.019 | 193.066 |
| | | $\bar{y}_{SK}$ | $a = 1, b = \beta_2(x)$ | 108.463 | 178.829 | 104.113 | 136.757 |
| | | $\bar{y}_{US1}$ | $a = \beta_2(x), b = C_x$ | 128.517 | 156.553 | 113.043 | 197.550 |
| | | $\bar{y}_{US2}$ | $a = C_x, b = \beta_2(x)$ | 113.065 | 199.197 | 106.257 | 149.504 |
| 1 | 0 | $\bar{y}_{ST}$ | $a = 1, b = \rho$ | 127.404 | 162.289 | 112.573 | 195.631 |
| 1 | $\frac{S_{YX}}{S_X^2}$ | $\bar{y}_{KC1}$ | $a = 1, b = 0$ | 481.283 | 57.570 | 514.286 | 125.884 |
| | | | $a = 1, b = C_x$ | 487.231 | 62.920 | 515.968 | 132.022 |
| | | | $a = 1, b = \beta_2(x)$ | 520.900 | 148.446 | 524.951 | 189.204 |
| | | | $a = \beta_2(x), b = C_x$ | 481.415 | 57.707 | 514.323 | 126.049 |
| | | | $a = C_x, b = \beta_2(x)$ | 514.286 | 123.659 | 523.258 | 177.845 |
| 1 | $\frac{S_{YX}}{S_X^2}$ | $\bar{y}_{reg}$ | -- | 526.316 | 208.264 | -- | -- |
| 0.610 | 70.689 | $\bar{y}_{GS}$ | -- | 863.816 | 214.473 | -- | -- |
| 0.610 | 70.689 | $\bar{y}_\lambda$ | -- | 863.816 | 214.473 | | |
| 0.971 | -0.540 | $\bar{y}_\lambda$ | -- | | | 863.816 | 214.473 |

**Conclusion**

In Table 2, we observe that most of the estimators considered in the paper performs better than usual ratio estimator ( except $\bar{y}_{KC1}$ which is inferior for some values). Also, we observe that many of these estimators are not efficient than regression estimator. The efficiency of the proposed estimator is same as that of Gupta and Shabbir [2] estimator and both the estimators performs better than regression estimator.

In this paper, we have suggested some estimators in simple random sampling without replacement using known values of some population parameters. It has been shown that many existing estimators are particular members of the proposed family of estimator. We studied the effect of various transformations of auxiliary information on the families of estimators.